\documentclass[12pt]{amsart}
\usepackage{geometry} 
\usepackage{graphicx}	
\usepackage{amsmath}								
\usepackage{amssymb}
\usepackage{mathabx}
\usepackage{stmaryrd}

\title[The Hybrid Discontinuous Galerkin method for elliptic problems]{The Hybrid Discontinuous Galerkin method for elliptic problems and applications
in vertical ocean-slice modeling}
\author[D. Azofeifa, M. A. Moreles, F. Velazquez]{Danalie Azofeifa, Miguel Angel Moreles, Federico Angel Velazquez-Mu\~{n}oz}

\address{D. Azofeifa, M. A. Moreles \\
              Centro de Investigaci\'{o}n en Matem\'{a}ticas\\ Jalisco s/n, Valenciana\\
Guanajuato, GTO 36240,Mexico\\
              \email{moreles@cimat.mx}    \\     
                         F. Velazquez \\
           Department of Physics, University of Guadalajara\\
Guadalajara 44100, Mexico
}
\date{} 

\begin{document}

\begin{abstract}
Classical numerical solutions of the Navier-Stokes equations applied to Coastal Ocean Modeling  are based on the Finite Volume Method  and the Finite Element Method. The Finite Volume Method guarantees local and global mass conservation. A property not satisfied by the Finite Volume Method. On the down side, the Finite Volume Method requires non trivial modifications to attain high order approximations unlike the Finite Volume Method. It has been contended that the Discontinuous Galerkin Method, locally conservative  and high order, is a natural progression for Coastal Ocean Modeling.  Consequently, as a primer we consider the vertical ocean-slice model with the inclusion of density effects. To solve these non steady Partial Differential Equations, we develop a pressure projection method for solution. We propose a Hybridized Discontinuous Galerkin solution for the required Poisson Problem in each time step. The purpose, is to reduce the computational cost of classical applications of the Discontinuous Galerkin method. The Hybridized Discontinuous Galerkin method is first presented as a general  elliptic problem solver.  It is shown that a high order implementation yields fast and accurate approximations on coarse meshes. 
\end{abstract}

\maketitle
\tableofcontents

\section{Introduction}
The numerical solution of the Navier-Stokes equations applied to Coastal Ocean Modeling (COM), is an active line of research. 
In Chen et al. \cite{Chen_etal2007}, a Finite Volume Method (FVM) for COM, referred as FVCOM is proposed and compared with finite difference (FD) models. Then an extension to  include non-hydrostatic effects is presented in Lai et al \cite{Lai_etal2010}.

 A motivation to develop FVCOM is the difficulty of dealing with irregular coastal geometries with FD, FVCOM applies naturally. Also, FVM guarantees local and global mass conservation. A property not satisfied by the very popular Finite Element Method (FEM). On the down side, 
the FVM requires non trivial modifications to attain high order approximations, unlike the FEM. Consequently, the next generation of COM is by means of methods that are locally conservative and high order. A natural choice is the Discontinuous Galerkin Method, Chen et al. \cite{Chen_etal2007}:
 \emph{Application of these methods to current finite element coastal ocean models could significantly improve computational accuracy and
efficiency as well as mass conservation}. These applications are under development, and there is an increasing literature on the subject. 
 See K\"{a}rn\"{a} et al \cite{Karnaetal} and Pan et al \cite{Panetal}. In the latter a non-hydrostatic extensions to a discontinuous finite element coastal ocean model is developed.

A drawback of DG methods is that in general it is more expensive that existing numerical methods, because DG methods  have many more (coupled) unknowns. A recent solution is to introduce hybridizable (hybridized) Discontinuous Galerkin Methods (HDG). The number of coupled unknowns is reduced, while retaining the attractive properties of the DG method, Tanh \cite{Bui_Thanh}.

Consequently, the purpose of our research, it to develop the HGD in the context of COM. Of particular interest are applications to non-hydrostatic modeling. As a primer, we consider the vertical ocean-slice model with the inclusion of density effects, K\"{a}mpf \cite{Kampf}, Lai et al. \cite{Lai_etal2010}. To solve these non steady PDE equations, we develop a pressure projection (PP) method for solution, see Almgren et al \cite{Almgren_etal}. It is well known that these PP methods require a Poisson solver in each time step, and we propose a HDG solution. Noteworthy, HDG methods do not require penalization unlike classical DG solutions for elliptic problems, see Riviere \cite{RiviereDG}.
Our implementations favor high order approximations, an attractive feature of DG methods. We do not aim for generality, comparisons with FVCOM are carried out only on rectangular meshes. It will become apparent that the application to more general meshes is straightforward.

The outline is as follows. In the Materials and Methods section, we develop an HDG scheme for the diffusion-advection-reaction equation. 
A distinctive step is to reduce the equation to a first order hyperbolic system.  We introduce an alternative to that of Bui-Thanh \cite{Bui_Thanh}. 
Also a PP method is developed for 2D vertical-slice modeling. We stress the use of the HDG solver of the underlying Poisson problem.
In Section 3 we consider pure elliptic problems in COM, namely tidal simulation in a semienclosed basin with tidal forcing at the open boundary.
Near resonance cases are considered. In Section 4 we address deep-water (short) surface gravity waves and density-driven currents, as vertical ocean-slice nonhydrostatic  models. Our findings and reflections on future work, are summarized in the Conclusions section.

\section{Materials and Methods}
\subsection{A HDG scheme for elliptic equations}
 Let us consider the diffusion-advection-reaction equation for the unknown function $p$,
\begin{equation}
-\nabla\cdot(\mathbf{K}\nabla p)+\mathbf{\beta}\cdot\nabla p+cp = f.
\label{D_A_R}
\end{equation}

It is assumed that the equation is of (strict) elliptic type, namely, the matrix function $\mathbf{K}$ is symmetric and positive definite in its domain of definition. 

Since the basic construction of the DG scheme is local, boundary conditions for well posedness shall be introduced in the corresponding step of the construction. 

\subsubsection{Function spaces}

Assume a domain $\Omega\subset\mathbb{R}^d$ partitioned into non-overlapping elements $I_j, \, j=1,\ldots, N$. Let $I$ be one of such elements. The boundary of $I$, $\partial I$ is the union of uniquely defined faces $e$ with outward pointing normal $\mathbf{n}$. The skeleton of the mesh, denoted by  $\mathcal{E}$ is the set of all faces. Two  subsets of interest are, the the set of boundary faces $\mathcal{E}^\partial$, and $\mathcal{E}^\circ=\mathcal{E}\setminus\mathcal{E}^\partial$ the set of interior faces.

Let $\mathbf{u},\mathbf{v}:\mathbb{R}^d\to\mathbb{R}^m$ smooth vectorial functions witn $m$ componentes. For $d-$dimensional set $G\subset\mathbb{R}^d$ we define the inner products
\[
(\mathbf{u},\mathbf{v})_G=\sum_{i=1}^m(\mathbf{u}_i,\mathbf{v}_i)_G\equiv\sum_{i=1}^m\int_G\mathbf{u}_i(\mathbf{x})\mathbf{v}_i(\mathbf{x})\, d\mathbf{x}.
\]
\[
\langle\mathbf{u},\mathbf{v}\rangle_{\partial G}=\sum_{i=1}^m\langle\mathbf{u}_i,\mathbf{v}_i\rangle_{\partial G}\equiv\sum_{i=1}^m\int_{\partial G}\mathbf{u}_i(\mathbf{x})\mathbf{v}_i(\mathbf{x})\, ds.
\]

For an element $I$, the latter reads
\[
\langle\mathbf{u},\mathbf{v}\rangle_{\partial I}=\sum_{e\subset \partial G, e\in\mathcal{E}}\int_e\mathbf{u}_i(\mathbf{x})\mathbf{v}_i(\mathbf{x})\, ds.
\]

Finally, we define $\mathcal{P}^p(I)$ the space of polynomials of degree less or equal than $p$ on $I$, and $\mathbf{V}^p_m(I)$ the cartesian product of $m-$copies of 
$\mathcal{P}^p(I)$. Similarly, $\mathcal{P}^p(e)$ and  $\mathbf{\Lambda}^p_m(e)$, for a face $e\in\mathcal{E}$.

\subsubsection{First order hyperbolic system}

Hereafter  we use freely the theory of hyperbolic partial differential equations, see Leveque  \cite{Randall_finite}. 

The first step is to reduce equation (\ref{D_A_R}), to a first order hyperbolic system.

Let us introduce a variant of the reduction in Bui-Thanh \cite{Bui_Thanh}.  We define the new variable (flux) $\mathbf{z}=-\mathbf{K}\nabla p$. 

For clarity of exposition, we consider two dimensional problems.

Since the matrix $\mathbf{K}$ is invertible, $\nabla p=-\mathbf{K}^{-1}\mathbf{z}$, we have the following system
\begin{eqnarray*}
\nabla p+\mathbf{K}^{-1}\mathbf{z}&=&\mathbf{0}\\
\nabla\cdot \mathbf{z}-\mathbf{\beta}\cdot \mathbf{K}^{-1} \mathbf{z}+cp&=&f
\end{eqnarray*}

Set $\mathbf{z}=(z^1, z^2), \mathbf{\beta}=(\beta_1, \beta_2)$, and define
\[
\textbf{u}:=\left(\begin{matrix}
z^1\\
z^2\\
p
\end{matrix}\right), \mathbf{A}_1:=\left(\begin{matrix}
0&0&1\\
0&0&0\\
1&0&0
\end{matrix}\right), \mathbf{A}_2:=\left(\begin{matrix}
0&0&0\\
0&0&1\\
0&1&0
\end{matrix}\right), 
\]
\[
\mathbf{B}:=\left(\begin{matrix}
k^{-1}_1&0&0\\
0&k^{-1}_2&0\\
-k^{-1}_1\beta_1&-k^{-1}_2\beta_2&c
\end{matrix}\right), \textbf{f}:=\left(\begin{matrix}
0\\
0\\
f
\end{matrix}\right).
\]
The system becomes
\begin{equation}
\frac{\partial}{\partial x}(\mathbf{A}_1\textbf{u})+\frac{\partial}{\partial y}(\mathbf{A}_2\textbf{u})+\mathbf{B}\textbf{u}=\textbf{f}.
\label{div_A}
\end{equation}

Let $\mathbf{n}$, be an arbitrary vector, and  $\mathbf{A}:=n_1\mathbf{A}_1+n_2\mathbf{A}_2$. It is readily seen that 
\[
\mathbf{A}=\left(\begin{matrix}
0&0&n_1\\
0&0&n_2\\
n_1&n_2&0
\end{matrix}\right)
\]
with real eigenvalues $\{0, \sqrt{n_1^2+n_2^2}, -\sqrt{n_1^2+n_2^2}\}$. Hence the matrix and the system are hyperbolic.\\

For later reference, let us write

\[
\mathbf{A}=\mathbf{R}\mathbf{D}\mathbf{R}^{-1}, \quad |\mathbf{A}|:=\mathbf{R}|\mathbf{D}|\mathbf{R}^{-1}
\]
where $\mathbf{R}=[\mathbf{r}_1, \mathbf{r}_2, \mathbf{r}_3]$ is the matrix of eigenvectors,  
\[
\mathbf{D}=diag(0, -\sqrt{n_1^2+n_2^2}, \sqrt{n_1^2+n_2^2})\equiv diag(\lambda_1, \lambda_2, \lambda_3),
\]
and  
\[
|\mathbf{D}|:=diag(|\lambda_1|, |\lambda_2|, |\lambda_3|).
\]

\bigskip

Finally, for any matrix $\mathbf{O}$, we denote by $\mathbf{O}_i,\, \mathbf{O}^j$ its row $i$, and column $j$ respectively.

\subsubsection{Local DG formulation with Godunov flux}

Let us define
\[
F_i(\textbf{u}):=((\mathbf{A}_1)_i\mathbf{u},(\mathbf{A}_2)_i\mathbf{u}), \quad i=1, 2,3,
\]
and 
\[
F(\textbf{u}):=\left[
\begin{array}{c}
F_1(\mathbf{u}) \\
F_2(\mathbf{u}) \\
F_3(\mathbf{u}) \\
\end{array}
\right].
\]

We apply operators component wise. For instance, for the divergence operator we have,
\[
\nabla\cdot F(\textbf{u}):=\left[
\begin{array}{c}
\nabla\cdot F_1(\mathbf{u}) \\
\nabla\cdot  F_2(\mathbf{u}) \\
\nabla\cdot  F_3(\mathbf{u}) \\
\end{array}
\right].
\]

We can write the system (\ref{div_A}) in the form
\begin{equation}
\nabla\cdot F(\mathbf{u})+\mathbf{B}\mathbf{u}=\mathbf{f},
\label{div_F}
\end{equation}

Compute the inner product  on $I$ of each side of (\ref{div_F}) with a test function $\mathbf{v}\in\mathbf{V}^p_3(I)$, to obtain
\begin{equation}\label{2_primera}
(\nabla\cdot F(\mathbf{u}),\mathbf{v})_I+(\mathbf{B}\mathbf{u},\mathbf{v})_I=(\mathbf{f},\mathbf{v})_I.
\end{equation}

Integrating by parts,
\begin{equation}\label{2_integrada}
-(F(\mathbf{u}),\nabla\mathbf{v})_I+\langle F(\mathbf{u})\cdot\mathbf{n},\mathbf{v}\rangle_{\partial I}
+(B\mathbf{u},\mathbf{v})_I=(\mathbf{f},\mathbf{v})_I,
\end{equation}

Continuity is not enforced at the boundary of adjacent elements. Therefore, the boundary term $F(\textbf{u})$ is replaced with a boundary numerical flux $F^*(u^-,u^+)$. As customary, the $-$ superscript denotes limits from the interior of $I$, and the $+$ superscript, limits from the exterior. In this context,  element $I$ is denoted by $I^-$ and the outer normal $\mathbf{n}$ by $\mathbf{n}^-$. 

A classical numerical flux is that of Godunov given by
\begin{equation}
F^*\cdot\mathbf{n}^-:=F(u^-)\cdot n^-+\vert \mathbf{A}\vert(u^--u^*).
\label{Gv_l}
\end{equation}
Here $u^*\equiv u^*(u^-,u^+)$ is the solution of a Riemann problem along the normal $\mathbf{n}^-$ of $I^-$.

The Godunov flux for the adjacent element $I^+$ on the same face of the boundary $\partial I^-$ is given by
\begin{equation}
F^*\cdot\mathbf{n}^+:=F(u^+)\cdot n^++\vert \mathbf{A}\vert(u^+-u^*).
\label{Gv_r}
\end{equation}
The following identity holds
\[
F^*\cdot\mathbf{n}^-+F^*\cdot\mathbf{n}^+=0.
\]
Or in terms of the  the jump operator 
\[
\llbracket (\cdot)\rrbracket = (\cdot)^-+(\cdot)^+,
\]
\[
  \llbracket   F(u)\cdot n+\vert \mathbf{A}\vert(u-u^*) \rrbracket  = 0.
\]

Combining adjacent fluxes
\[
F^*\cdot\mathbf{n}^-=\frac{1}{2}\left[F(u^-)+F(u^+)\right]\cdot n^-+\frac{1}{2}\vert \mathbf{A}\vert(u^--u^+).
\]
This is the symmetric form of the Godunov flux used in upwind DG. It couples the unknowns of the adjacent elements, and hence the unknowns of all elements.

\subsubsection{Hybrid flux}
Observe that the upwind fluxes (\ref{Gv_l}), (\ref{Gv_r})  depend on the DG unknowns of only one side of a face and the single-valued solution $u^*$ of the Riemann problem. If $u^*$ is given, the numerical flux is completely determined using only information from either side of the face. Moreover,  we then can solve for $u$ element-by-element independent of each other.

To hybridized the flux, and break the coupling, $u^*$ is regarded as an extra unknown to be solved on the skeleton of the mesh instead of using the Riemann/upwind state which couples the local unknown $u$. Renaming $u^*$ as $\hat{u}$ and $F^*$ as $\hat{F}$, we are led to
\begin{equation}
\hat{F}\cdot\mathbf{n}:=F(u)\cdot n+\vert \mathbf{A}\vert(u-\hat{u}).
\label{HDG_flux}
\end{equation}
This is  the hybridized upwind flux or HDG flux.

In summary, for each element $I$, the DG local unknown $u$ and the extra \emph{trace} unknown $\hat{u}$ need to satisfy
\begin{equation}
-(F(\textbf{u}),\nabla \textbf{v})_I+\langle \hat{F}(\mathbf{u}^-,\hat{u})\cdot\mathbf{n}^-,\mathbf{v}\rangle_{\partial I}+(\mathbf{B}\textbf{u},\mathbf{v})_I=(\textbf{f},\mathbf{v})_I,
\quad \mathbf{v}\in\mathbf{V}^p_3(I).
\label{interior_eq}
\end{equation}

This is complemented with a weak jump condition in the skeleton. Namely, for all $e\in\mathcal{E}$,
\begin{equation}
 \langle \llbracket   \hat{F}\cdot\mathbf{n} \rrbracket,\mathbf{w}\rangle_e  = 0, \quad \mathbf{w}\in \mathbf{\Lambda}^p_m(e).
 \label{bd_eq}
\end{equation}

\subsubsection{The discrete problem}

Let us solve equation (\ref{interior_eq}) for the $\mathbf{u}$ terms,  and  (\ref{bd_eq})
for the $\hat{\mathbf{u}}$ terms. We obtain,
\begin{equation}
-(F(\textbf{u}),\nabla \textbf{v})_I+\langle F(\mathbf{u}^-)\cdot\mathbf{n}^-+\vert \mathbf{A}\vert\mathbf{u}^-,\mathbf{v}\rangle_{\partial I}+(\mathbf{B}\textbf{u},\mathbf{v})_I=(\textbf{f},\mathbf{v})_I+
\langle\vert \mathbf{A}\vert\hat{\mathbf{u}},\mathbf{v}\rangle_{\partial I}.
\label{interior_eq_u}
\end{equation}

\begin{equation}
 \langle   2\vert \mathbf{A}\vert\hat{\mathbf{u}},\mathbf{w}\rangle_e =
 \langle \llbracket   F(\mathbf{u} )\cdot\mathbf{n} +\vert \mathbf{A}\vert\mathbf{u} \rrbracket ,\mathbf{w}\rangle_e.
 \label{bd_eq_hat_u}
\end{equation}

Let $\lbrace \mathbf{N}_j: j=1,\ldots, P_I\rbrace$ be a basis of  $\mathbf{V}^p_3(I)$. Hence
\[
\mathbf{u} = \sum_{j=1}^{P_I}\mathbf{u}_j\mathbf{N}_j
\]

Since $F(\mathbf{u} )=\mathbf{A}_1\mathbf{u}+\mathbf{A}_2\mathbf{u}$ and $A= (n_I)_1A_1+(n_I)_2A_2$,  
for $ i=1,\ldots, P_I$, equation (\ref{interior_eq_u})  reads
\begin{equation}
\begin{array}{c} 
\sum_{j=1}^{P_I}\left[
(-(\mathbf{A}_1+\mathbf{A}_2)\mathbf{N}_j,\nabla\cdot\mathbf{N}_i)_I+
(B\mathbf{N}_j,\mathbf{N}_i)_I+
\langle (\mathbf{A}+\vert \mathbf{A}\vert)\mathbf{N}^-_j,\mathbf{N}^-_i\rangle_{\partial I}
\right]
\mathbf{u}_j \\ = \\
(\textbf{f},\mathbf{N}_i)_I+
\langle\vert \mathbf{A}\vert\hat{\mathbf{u}},\mathbf{N}^-_i\rangle_{\partial I}
\end{array}
\label{local_I}
\end{equation}

Similarly, let $\lbrace \mathbf{M}_j: j=1,\ldots, Q_e\rbrace$ be a basis of  $\mathbf{\Lambda}^p_3(e)$.
For $ i=1,\ldots, Q_e$, equation (\ref{bd_eq_hat_u})  reads
\begin{equation}
\sum_{j=1}^{Q_e} \langle   2\vert \mathbf{A}\vert\mathbf{M}_j,\mathbf{M}_i\rangle_e\hat{\mathbf{u}}_j =
 \langle \llbracket   F(\mathbf{u} )\cdot\mathbf{n} +\vert \mathbf{A}\vert\mathbf{u} \rrbracket ,\mathbf{w}\rangle_e.
 \label{local_e}
\end{equation}

Solving for $\mathbf{u}$ on each element form (\ref{local_I}) and substituting in the corresponding edges in equation (\ref{local_e}), we are led to a sparse linear system for the hatted unknown on the skeleton.  

Then one solves the small and independent linear systems (\ref{local_I}) for the 
local variables on each element. The latter can be done in parallel.

\subsubsection{Element basis functions}

The DG method is $H^p$ adaptative, it is straightforward to make the order of approximation element dependent. The method is also suitable for irregular geometries and unstructured meshes. Nevertheless, we do not aim for generality. We shall test on benchmark problems in rectangular geometries and regular meshes.

For one dimensional problems we use second order nodal functions on the reference interval $[-1,1]$,
\[
N_0(\xi)=\frac{\xi^2-\xi}{2},\quad N_1(\xi)=1-\xi^2,\quad N_2(\xi)=\frac{\xi^2+\xi}{2}.
\]

An easy 2D extension is achieved by considering the basis functions
\[
N_{ij}(\xi,\zeta)=N_i(\xi)\, N_j(\zeta), \quad i,j=0,1,,2,
\]
defined on the reference square $[-1,1]\times[-1,1]$.

This basis will help to illustrate the performance of a high order DG method.

\subsection{A pressure projection method with HDG Poisson solver}

Pressure projection methods are well known, a thorough study is presented in Almgren et al \cite{Almgren_etal}. As a primer, we introduce a PP method in a pressure splitting framework to be applied in 2D vertical-slice modeling. Our purpose is to stress the use of the HDG solver of the underlying Poisson problem.

Let us consider unsteady, constant density, incompressible flow. Let  $\rho_0$ be the reference density. The flow is governed by the momentum and continuity equations
\[
\frac{\partial\mathbf{u}}{\partial t}=-\frac{1}{\rho_0}\nabla P-(\mathbf{u}\cdot\nabla)\mathbf{u},
\]
\[
\nabla\cdot\mathbf{u}=0.
\]

 Let $\Delta t$ be the time step. Assume the velocity $\mathbf{u}^n$ and pressure $P^n$ are given at time $t_n$. 

Split pressure $P^{n+1}$ in the form
\[
P^{n+1}=p^n+q^{n+1}
\]
where $p^n$, $q^{n+1}$  are hydrostatic and non hydrostatic pressure respectively. The latter is of
the form
\[
q^{n+1}=q^n+\delta q^{n+1}.
\]

The non hydrostatic correction $\delta q^{n+1}$ is to be determined. 

The main steps are:

\bigskip

\noindent\textbf{Step 1.} Construct an intermediate velocity $\mathbf{u}^{n+1/2}$ by advancing the momentum equations,
\[
\mathbf{u}^{n+1/2}=\mathbf{u}^{n}-\frac{\Delta t}{\rho_0}\nabla P^n-(\mathbf{u}^n\cdot\nabla)\mathbf{u}^n
\]

\noindent\textbf{Step 2.} Pressure correction. For suitable boundary conditions, HDG solve the Poisson equation
\[
-\Delta (\delta q^{n+1})=-\frac{\rho_0}{\Delta t}\nabla\cdot\mathbf{u}^{n+1/2}.
\]

\noindent\textbf{Step 3.} Update the divergence free velocity field,
\[
\mathbf{u}^{n+1}=\mathbf{u}^{n+1/2}-\frac{\rho_0}{\Delta t}\nabla(\delta q^{n+1}).
\]

\bigskip

It is straightforward to modify this scheme for more general unsteady equations. We show some examples below.

\section{Tidal simulation in semienclosed basin}

 We consider two benchmark problems that lead to elliptic equations. A comparison is made with the solution of  FVCOM as presented in Chen et al \cite{Chen_etal2007}.  Therein, FVCOM is applied for modeling of tidal simulation in semienclosed basin with tidal forcing at the open boundary under nonresonance and near resonance conditions.  Here we apply HDG to illustrate the accurate simulation of the troublesome near resonance case. 
 
 \subsection{A Rectangular Channel}

Consider a fluid layer of uniform density that propagates  along a channel aligned with the $x-$direction. More precisely, a semienclosed narrow channel with length $L$ and variable depth $H(x)$ and a closed boundary at $x=L_1$ and an Open Boundary at $x=L$.

\begin{figure}[h]
\center
\includegraphics[scale=0.6]{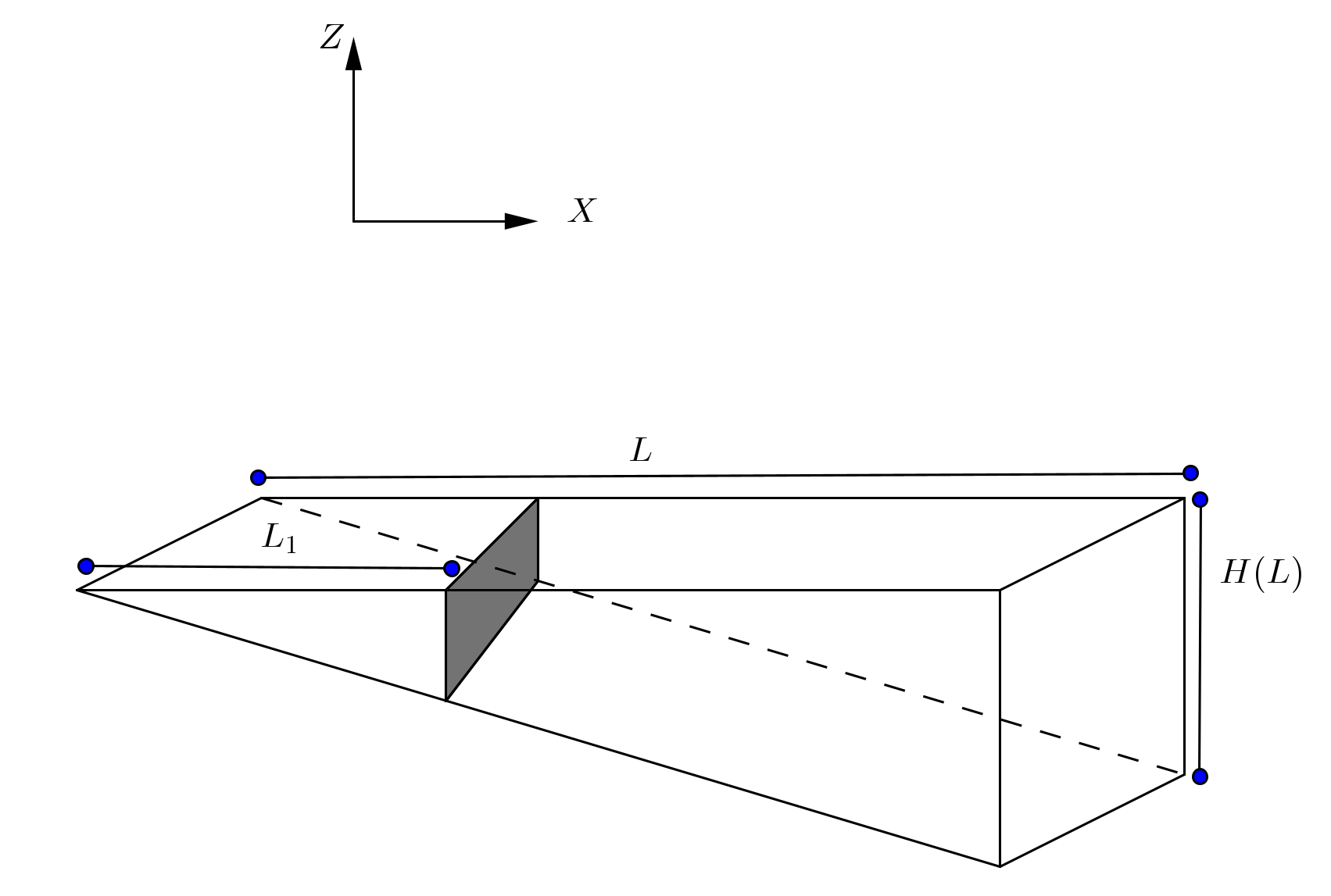}\\
\caption{\textit{Configuration of the semienclosed channel.}}
\label{canal_1d}
\end{figure}

Neglecting Coriolis force and advection of momentum,  the governing equations modeling tidal waves propagation in the semienclosed channel (see Figure \ref{canal_1d}) are given as
\[
\frac{\partial u}{\partial t}+g\frac{\partial \zeta}{\partial x}=0;\quad \frac{\partial \zeta}{\partial t}+g\frac{\partial uH}{\partial x}=0;
\quad (x,t)\in (a,b)\times(0,T).
\]
Here,  $H(x)$ is the total water depth, $g$ is acceleration due to gravity, $u$ is speed in the $x-$direction, and $\zeta$ is sea-level elevation .\\

Assuming harmonic solutions,
\[
\zeta=\zeta_0(x)e^{-i\sigma t},\quad u=u_0(x)e^{-i\sigma t},
\]
we obtain the ordinary second order equation for $\zeta_0$ 
\begin{equation}
(H\zeta_0')'+\frac{\sigma^2}{g}\zeta_0=0.
\label{Eq_1D}
\end{equation}

We specify a periodic tidal forcing with amplitude $A$ at the mouth of the channel, 
\[
\zeta_0(L)=A,
\]
and a no-flux boundary condition at the wall,
\[
H(L_1)\zeta_0'(L_1)=0.
\]

Let the water depth decrease linearly toward the end of the channel, so that $H(x)$ can be written as
$$H(x)=\frac{xH(L)}{L}.$$

It is readily seen that (\ref{Eq_1D}) is a Bessel's equation. With the given boundary conditions,
the analytic solution is given by 

$$\zeta_0(x)=A\frac{Y_0^{\prime}(2k\sqrt{L_1})J_0(2k\sqrt{x})-J_{0}^{\prime}(2k\sqrt{L_1})Y_{0}(2k\sqrt{x})}{Y_0^{\prime}(2k\sqrt{L_1})J_{0}(2k\sqrt{L})-J_0^{\prime}(2k\sqrt{L_1})Y_0(2k\sqrt{L})}$$
where 
$$k:=\frac{\sigma\sqrt{L}}{\sqrt{gH(L)}},$$
$J_0, Y_0$ are the Bessel's functions the degree zero and one respectively.\\

To compare with the HDG approximation, we solve the  hyperbolic system,
\begin{eqnarray*}
\zeta_0^{\prime}+H^{-1}z&=&0\\
z^{\prime}-\frac{\sigma^2}{g}\zeta_0&=&0.
\end{eqnarray*} 
Here  $z=-H\zeta_0^{\prime}$, (\ref{Eq_1D}).

The following parameters are considered for a channel very close to resonance. 

\[
L=300km,\, L_1=10km, \,H(L)=20.1m,\, \sigma=\frac{2\pi}{12.42\cdot 3600s}, \,A=1cm.
\]

The HDG method is applied using nodal polynomials of degree 2 ($\mbox{HDG}_2$). It is compared with the FV method and the analytic solution. 
For consistency with the FV method, we consider the approximation at the middle point $x_i$ of the element  $I_i=[x_{i-\frac{1}{2}}, x_{i+\frac{1}{2}}]$. 
We compute the relative root mean square error  for $\zeta_0^R$ the analytic solution and $\zeta_0$ the numerical solution. Namely,
\[
\mathcal{E}_2:=\frac{||\zeta_0-\zeta_0^R||_2}{||\zeta_0^R||_2}.
\]

As illustrated in Table \ref{tabla_1}, a second order approximation with HDG in a coarse resolution, is of greater quality than FV. 

\begin{table}
\centering
\begin{tabular}{|c|c|c|}
\hline
Elements & $\mathcal{E}_2(\mbox{HDG}_2)$ & $\mathcal{E}_2(FV)$\\
\hline
10 & 0.180784 & 22.3052\\
\hline
20 & 0.0748017 & 0.585413\\
\hline
40 & \textbf{0.0235224} & 0.553751\\
\hline
80 & 0.00684013 & 0.43336\\
\hline
160 & 0.00187985 & 0.297118\\
\hline
320 & 0.000495726 & 0.181964\\
\hline
640 & 0.000127432 & 0.102475\\
\hline
1280 & $0.0000324$ & \textbf{0.0546902}\\
\hline
\end{tabular}
\label{tabla_1}
\caption{Relative root mean square error for $\mbox{HDG}_2$ and FV.}
\end{table}

\begin{figure}[h]
\centering
\includegraphics[scale=0.6]{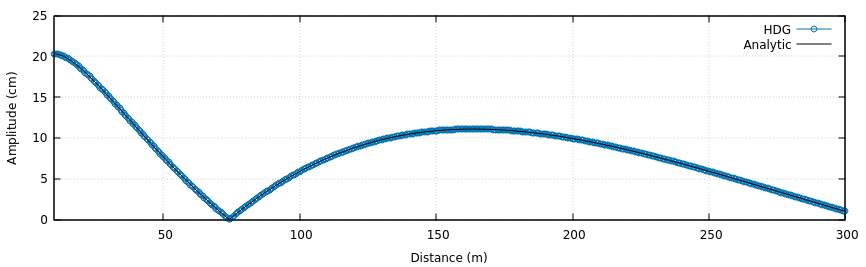}
\caption{Solution for the rectangular channel near resonance case, with HDG nodal with polynomials of degree 2 and 80 elements.}
\label{compara_HDG_finite_volume}
\end{figure}

The analytic solution describes a standing wave with a node point  near the closed side of the channel. As expected, the reproduction of  these features by HDG is accurate. See Figure \ref{compara_HDG_finite_volume}

\bigskip

\noindent\textbf{Remark. }As pointed out in Chen et al \cite{Chen_etal2007}, regardless of the numerical method, a proper selection of horizontal resolution recover accurately this tidal resonance problem.  We argue that the accuracy attained  by HDG in coarse meshes yields a better choice.

\subsection{A Sector Channel}

Now consider a flat bottom channel in the form of a semicircular section, which in polar coordinates is defined from $0$ to $L$ in the radial direction and from $\alpha/2$ to $\alpha/2$ in the angular direction.

\begin{figure}[h]
\centering
\includegraphics[scale=0.35]{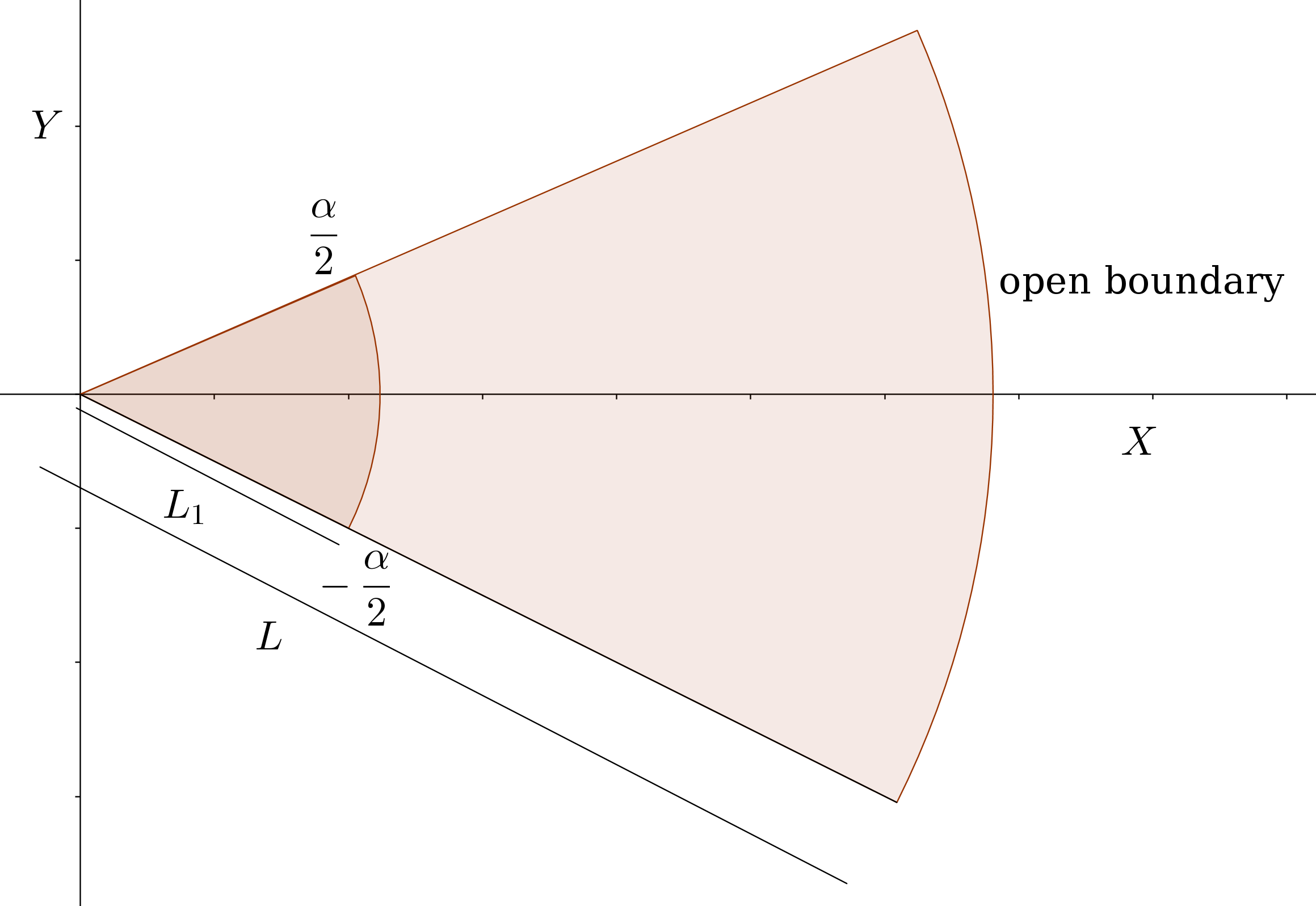}
\caption{Semienclosed sector channel in the polar region $L_1\leq r\leq L$,  $-\frac{\alpha}{2}\leq\theta\leq \frac{\alpha}{2}$. The sector is open at $r=L$ and closed elsewhere.}
\label{channel_2D}
\end{figure}

The semicircular line of radius $L$ corresponds to an open border, while along the semicircular line of radius $L_1$ and the two sides, they are closed, (Figure \ref{channel_2D}). 

The following equations govern the non rotating tidal oscillation,
\begin{eqnarray}
\frac{\partial V_r}{\partial t}&=&-g\frac{\partial\eta}{\partial r},\label{canal_2d_e1}\\
\frac{\partial V_{\theta}}{\partial t}&=&-g\frac{\partial\eta}{r\partial\theta},\\
\frac{\partial\eta}{\partial t}+\frac{\partial r V_rH_0}{r\partial r}+\frac{\partial V_0H_0}{r\partial\theta}&=&0.\label{canal_2d_e2}
\end{eqnarray}

$H_0$ is the constant water depth, $V_r, V_{\theta}$ are the radial and angular $r, \theta$ velocity components, and $\eta$ is the free surface water elevation.\\

Assuming a harmonic solution,

$$\eta=Re(\eta_0(r, \theta)e^{-it(\omega t-\frac{\pi}{2})}),$$
we can reduce the equations (\ref{canal_2d_e1}-\ref{canal_2d_e2}) to an elliptic equation
\begin{equation}\label{equa_canal_2d}
\frac{\partial^2\eta_{0}}{\partial r^2}+\frac{1}{r}\frac{\partial\eta_0}{\partial r}+\frac{1}{r^2}\frac{\partial^2\eta_0}{\partial\theta^2}+\frac{\omega^2}{gH_0}\eta_0=0.
\end{equation}

The physical boundary conditions are as follows
\begin{enumerate}
\item At the open mouth of the channel, a harmonic tidal forcing is assumed,
\[
\eta_0(r, \theta)=\overline{A}\cos\left(m\pi\frac{\theta+\frac{\alpha}{2}}{\alpha}\right),\quad \lbrace L\rbrace\times\left(-\frac{\alpha}{2},\frac{\alpha}{2}\right),
\] 
\item On the solid walls, null flux is prescribed,
\[
-K\nabla\eta_0(r, \theta)=0,\quad \lbrace L_1\rbrace\times\left(\frac{-\alpha}{2},\frac{\alpha}{2}\right)\bigcup
 (L_1,L)\times\left\lbrace-\frac{\alpha}{2}\right\rbrace\bigcup
 (L_1,L)\times\left\lbrace\frac{\alpha}{2}\right\rbrace.
\]
\end{enumerate}

The analytic solution of this boundary-value problem is:
$$\eta_0(r, \theta)=\overline{A}\frac{Y^{\prime}_v(L_1\kappa)J_v(r\kappa)-J^{\prime}_v(L_1\kappa)Y_v(r\kappa)}{Y^{\prime}_v(L_1\kappa)J_v(L\kappa)-J^{\prime}_v(L_1\kappa)Y_v(L\kappa)}\cos\left(\frac{m\pi}{\alpha}(\theta+\frac{\alpha}{2})\right)$$ 
where
$$v=\frac{m\pi}{\alpha}, \kappa=\frac{\omega}{\sqrt{gH_0}},$$
$J_v, Y_v$ are respectively the $v$ th-order Bessel function of the first and second type.

Let us show the HDG solution in the rectangular domain $(L_1,L]\times (-\alpha/2,\alpha/2)$. \\

Let $\nabla$ denote the gradient with respect to $(r, \theta)$ and let
$$K:=\left(\begin{matrix}
1&0\\
0&\frac{1}{r^2}
\end{matrix}\right).$$

Equation (\ref{equa_canal_2d}) becomes the  diffusion-advection-reaction equation,
\begin{equation}
-\nabla\cdot(K\nabla\eta_0)-(r^{-1}, 0)\cdot\nabla\eta_0-\frac{\omega^2}{gH_0}\eta_0=0.
\end{equation}

We apply the fourth order ($\mbox{HDG}_4$) scheme developed above for a near resonance case.  The geometric parameter values are:
\[
H_0=1 m,\, \alpha=\frac{\pi}{4}, \,L_1=90 km, \, m=1.0, \, \omega=\frac{2\pi}{12.42\cdot 3600 }\frac{1}{s}, \, L=158 km, \overline{A}=1 cm.
\]

\begin{figure}
\centering
a)\hspace{7cm}b)\\
\includegraphics[scale=0.45]{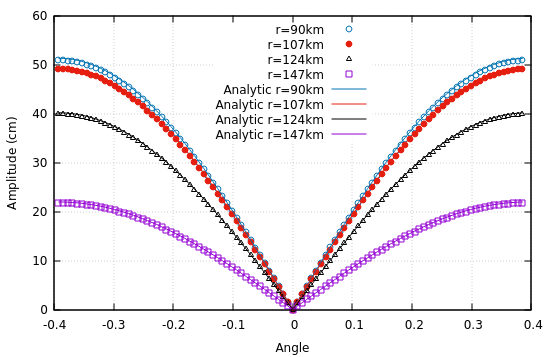}
\includegraphics[scale=0.45]{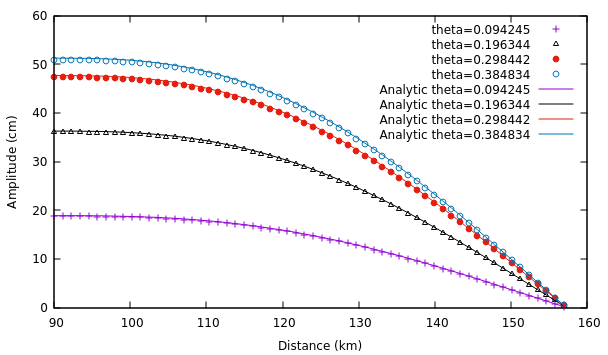}
\caption{Comparison between the analytic solution to the sector problem with the  $\mbox{HDG}_4$ solution using $30\times 50$ elements. a) Solution for some  $r$ values, b) Solution for some $\theta$ values.}
\label{figure_canal2d_2}
\end{figure}

In Figure (\ref{figure_canal2d_2}) the analytic and numerical solution are compared.  The relative mean square error is shown in Table \ref{tabla2_HDG4_HDG2_FV_1}. 

\begin{table}[h]
\centering
\begin{tabular}{|c|c|c|}
\hline
Elements & $\mathcal{E}_2(\mbox{HDG}_{4})$ & $\mathcal{E}_2(FV)$\\
\hline
$5\times 5$ & 0.441213 & 0.920939\\
\hline
$10\times 10 $ & \textbf{0.009879} & 0.887757\\
\hline
$35\times 35$ & 0.0026408 & 0.803039\\
\hline
$40\times 40$ & 0.00115041 & \textbf{0.794262}\\
\hline
\end{tabular}
\caption{Relative root mean square error for $\mbox{HDG}_4$ and FV.}
\label{tabla2_HDG4_HDG2_FV_1}
\end{table}

\bigskip

\noindent\textbf{Remark. }Noteworthy, the $10\times 10$ $\mbox{HDG}_4$ solution is of grater quality that the $40\times 40$ FV solution.
On regards to execution time, the former is attained in half a second in a personal laptop. A solution of the same quality would require in the order of minutes with the Finite Volume Method in a much finer mesh.

\section{Vertical ocean-slice nonhydrostatic  models}

Our proposal is to use a high order HDG Poisson solver to accelerate pressure projection numerical methods. As a primer, we solve two simple unsteady problems in Coastal Ocean Modeling. 

Let us consider the ocean as a vertical slice. As customary, flow and gradients of variables normal to this plane are assumed to vanish, and
the Coriolis force is ignored. We address deep-water (short) surface gravity waves and density-driven currents.

Modeling is considered in the Cartesian coordinate system, in which $x, y$ denote the horizontal coordinates, the slice is along the $x-$axis. Changes are constant on $y-$cross sections, all the normal forces to this slice are neglected. The vertical axis $z$ points upward to the undisturbed surface water located at $z=0$, $\eta$ is the free surface water, $h_0$ is the undisturbed water depth, and $h$ the total depth.  

Assuming a flat bottom, the time space domain of interest is $\Omega_T=[0, T_f]\times[L_a, L_b]\times [-h_0, 0]$. 
For simulation, a uniform rectangular mesh is considered with $\Delta x\times\Delta z$ elements.


\subsection{Surface Standing Waves in a Deep Basin}

The benchmark that follows is solved with a non hydrostatic version of FVCOM in Lai et al \cite{Lai_etal2010}. 
It was selected to test  the stability and accuracy of the split mode explicit non hydrostatic
algorithm and the non hydrostatic pressure Poisson solver.  The surface boundary condition requires an approximate treatment.
 Our pressure projection method with HDG Poisson solver is developed as an alternative.

 Consider the governing equations for a standing wave under linear, inviscid, non rotating conditions,
\begin{eqnarray*}
\frac{\partial u}{\partial t}&=&\frac{-1}{\rho_0}\frac{\partial P}{\partial x},\\
\frac{\partial w}{\partial t}&=&\frac{-1}{\rho_0}\frac{\partial P}{\partial z},\\
\frac{\partial u}{\partial x}+\frac{\partial w}{\partial z}&=&0.
\end{eqnarray*} 
Here, $u, w$ are the velocities in direction $x$ and $z$ respectively, $P$ is the dynamic pressure 
\[
P=p+q,
\]
$p$ the hydrostatic and $q$ the non hydrostatic pressure. $\rho_0$ is the mean density. 
The latter equation simplifies to $P = q$.

The free surface water, $\eta$, is given by the volume-conservation,
\begin{equation}\label{prognostic}
\frac{\partial\eta}{\partial t}=-\frac{\partial}{\partial x}(h<u>);\quad  <u>(x):=\frac{1}{h}\int_{z_0}^{z_0+h}u(t, x, z)dz.
\end{equation}

To relate the sea level elevation with the non hydrostatic dynamic pressure, the hydrostatic approximation is used,
\begin{equation}\label{prognostic_presion}
q_s=\rho_0g\eta\Longrightarrow \frac{\partial q_s}{\partial t}=-\rho_0 g\frac{\partial}{\partial x}(h<u>).
\end{equation}
Here, $g:=9.81 \frac{m}{s^2}.$ \\  

Hence, the full set of equations are,
\begin{eqnarray}
\frac{\partial u}{\partial t}&=&\frac{-1}{\rho_0}\frac{\partial q}{\partial x},\label{NO_reducido_1}\\
\frac{\partial w}{\partial t}&=&\frac{-1}{\rho_0}\frac{\partial q}{\partial z},\label{NO_reducido_1_1}\\
\frac{\partial u}{\partial x}+\frac{\partial w}{\partial z}&=&0,\label{NO_reducido_1_2}\\
\frac{\partial q_s}{\partial t}&=&-\rho_0 g\frac{\partial}{\partial x}(h<u>).\label{NO_reducido_2}
\end{eqnarray} 

\subsubsection{Pressure splitting.} Here we apply the PP method to the velocity field
\[
\mathbf{u}=\left(
\begin{array}{c}
u \\ w
\end{array}
\right).
\]

It is left to update the  correction pressure at the surface $\delta q_s^{n+1}$. 

The approximation is made using the fact that
$$q^{n+1}_s=q^{n}_s+\delta p^{n+1}_s\Longrightarrow \delta q^{n+1}_s=q^{n+1}_s-q^{n}_s.$$
It follows at once that
\begin{eqnarray*}
\delta q^{n+1}_s&=&\triangle t\frac{q^{n+1}_s-q^{n}_s}{\triangle t}\\
&\approx&\triangle t\frac{\partial q_s^{n}}{\partial t}\\
&=&-\triangle t\rho_0 g\frac{\partial}{\partial x}(h<u^n>).
\end{eqnarray*}

Approximating the derivative with respect to $x$ we have
\begin{equation}\label{correccion_presion_simple}
\delta q^{n+1}_s=-\frac{\Delta t}{\Delta x}\rho_0 g\left[\int_{-h_0}^{\eta^n}u(t^{n}, x_{k}, z)dz-\int_{-h_0}^{\eta^n}u(t^{n}, x_{k-1}, z)dz\right].\\
\end{equation}

The following numerical results was solved the HDG with nodal polynomials of degree $\vartheta=2$. 

\subsubsection{A test example}
Lai et al. \cite{Lai_etal2010}  run FVCOM in a closed rectangular channel  $L=10$ meters long and $H=10$ meters deep. 

The initial condition are,
\begin{eqnarray*}
u(0, x, z)&=&0,\\
w(0, x, z)&=&0,\\
q(0, x, z)&=&\rho_0 g\eta_0\frac{\cosh(\kappa(-z+H))}{\cosh(\kappa H)}\cos(\kappa x).\\
\end{eqnarray*}

For the next time steps, we set for the Poisson equation
\begin{eqnarray*}
q(t>0, x, 0)&=&\rho_0 g\frac{\cosh(\kappa(H))}{\cosh(\kappa H)}\eta(x),\\
\nabla q(t>0, x, z>0)\cdot\overrightarrow{\textbf{n}}&=0&,
\end{eqnarray*}
where $\rho_0=1000\, Kg/m^3$. 
A unimodal standing wave of free surface perturbation is considered, that is, 
\[
\eta(x)=\eta_0\cos(\kappa x)\cos(\omega t).
\]

The analytical solution of this problem can be found for instance in Jankowski \cite{Jankowski}. It is a non hydrostatic deep water wave with phase speed
$\omega=\sqrt{g\kappa\tanh(\kappa H)}$.

In Lai et al. \cite{Lai_etal2010} a non overlapping unstructured triangular grid is created  by dividing each square into two
triangles, with a triangle's side length of 0.25 m. Consequently we work with a square mesh with $dx = dz = 0.25$. 
In our simulation $\kappa=\frac{\pi}{L}=\frac{\pi}{10}$, $\eta_0 = 0.1$ small enough so that $\frac{\eta_0}{H}<< 1$. 

As in Lai et al. \cite{Lai_etal2010} to test the stability of the program, it is run for 10 minutes. Therein the time step is $0.05s$, the high order used in HDG allows a larger step, $dt = 0.5s$.  

\begin{figure}[h]
\centering
\includegraphics[scale=0.7]{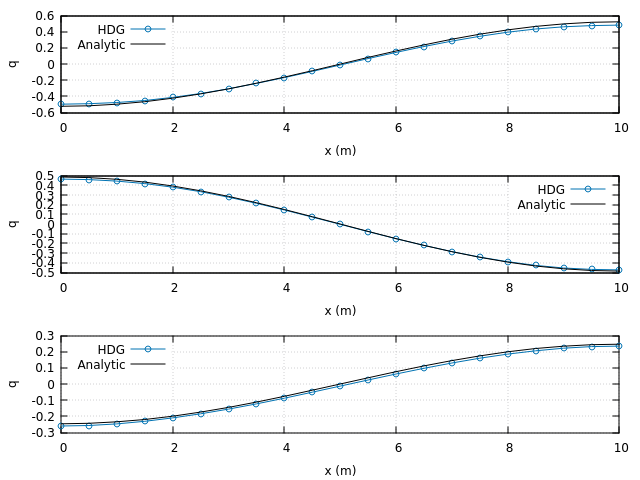}\\
\caption{Comparison between the HDG and analytic solutions  of  pressure at a depth of  $z=1$ meters for $t=2, 4, 10$ minutes respectively.}
\label{basin_images_corte}
\end{figure}

In Figure (\ref{basin_images_corte}) we show some snapshots of  pressure, scaled to $[-1, 1]$. 

\subsection{Density driven flow}
Now we consider density effects in the vertical ocean-slice model by adding an advection-diffusion equation for density and  the reduced-gravity force in the
vertical momentum equation. The governing equations can be written as: 
\begin{eqnarray}
\frac{\partial u}{\partial t}&=&-\frac{1}{\rho_0}\frac{\partial P}{\partial x}-(u, w)\cdot\nabla u,\label{e_density_i}\\
\frac{\partial w}{\partial t}&=&-\frac{1}{\rho_0}\frac{\partial P}{\partial z}-\frac{\rho-\rho_0}{\rho_0}g-(u, w)\cdot\nabla w,\\
\frac{\partial\rho}{\partial t}&=&\nabla\cdot\left(\left(\begin{matrix}
K_h&0\\
0&K_z
\end{matrix}\right)\nabla\rho\right)-(u, w)\cdot\nabla\rho,\\
\nabla\cdot(u, w)&=&0,\\
P_s&=&-\rho_0g\frac{\partial}{\partial x}(h<u>).\label{e_density_f}
\end{eqnarray}
Splitting the pressure in to hydrostatic and non hydrostatic pressure respectively,
$$P=p+q,$$
and using the hydrostatic approximation
$$-\frac{1}{\rho_0}\frac{\partial p}{\partial z}-\frac{\rho-\rho_0}{\rho_0}g=0$$
The set of equations (\ref{e_density_i}-\ref{e_density_f}) can be written as
\begin{eqnarray}
\frac{\partial u}{\partial t}&=&-\frac{1}{\rho_0}\frac{\partial (p+q)}{\partial x}-(u, w)\cdot\nabla u,\\
\frac{\partial w}{\partial t}&=&-\frac{1}{\rho_0}\frac{\partial q}{\partial z}-(u, w)\cdot\nabla w,\\
\frac{\partial\rho}{\partial t}&=&\nabla\cdot\left(\left(\begin{matrix}
K_h&0\\
0&K_z
\end{matrix}\right)\nabla\rho\right)-(u, w)\cdot\nabla\rho,\\
\nabla\cdot(u, w)&=&0,\\
P_s&=&-\rho_0g\frac{\partial}{\partial x}(h<u>),\\
\frac{\partial p}{\partial z}&=&-(\overline{\rho}-\rho_0)g, p(z=0)=0.
\end{eqnarray}

Vertical elements are indexed by $i$. Thus, we have $\overline{\rho}:=0.5(\rho_{i-1}+\rho_{i})$ in element $i$.

\subsubsection{Pressure splitting.}
To solve the last set of equation, we apply again pressure splitting. To complete the scheme, water's density is updated by a simple finite difference,
\[
\rho^{n+1}=\rho^n+\Delta t\nabla\cdot\left(\left(\begin{matrix}
K_h&0\\
0&K_z
\end{matrix}\right)\nabla\rho^{n}\right)-\Delta t(u^{n}, w^{n})\cdot\nabla\rho^{n}.
\]
whereas the hydrostatic pressure
\[
\frac{\partial p}{\partial z}=-(\overline{\rho}-\rho_0)g, p(z=0)=0,  \overline{p}:=0.5(p_{i-1}+p_{i})
\]

\subsubsection{A density driven flow with variable bottom topography}

\begin{figure}[h]
\centering
\includegraphics[scale=0.3]{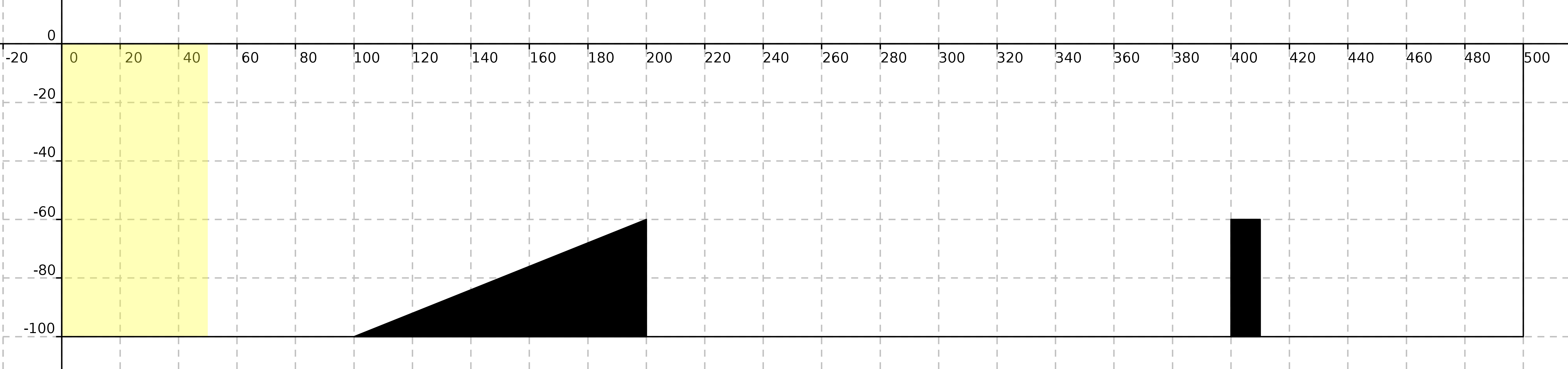}
\caption{Bathymetry for the variable density problem. The left rectangular (yellow color) is a column of water with density $1029 kg/m^3$, the triangle and the right rectangular (both in black color) are part of the sea floor. The remain domain is water with density $1028kg/m^3$.}
\label{del_densidad_variable_batimetria}
\end{figure}

As a final example, let us consider a closed channel problem solved by Finite Differences in K\"{a}mpf \cite{Kampf}. 
The closed channel is initially composed of two vertical layers of water with constant but distinct density.  The model is forced via prescription of a layer of dense water that initially leans against the left boundary. The variable bottom topography includes a ramp and a vertical bar, see Figure \ref{del_densidad_variable_batimetria}.

The closed channel occupies the rectangle  $\Omega=[0, 500]\times[-100,0]$. The reference density is given by $\rho_0=1028 kg/m^3$, and a column of water with density $\rho_1=1029 kg/m^3$  is included in the region $\Omega_1=[0, 50]\times[-100, 0]$.  The system starts at rest.

For the horizontal and vertical density diffusivities we set $K_h=K_z=10^{-4}m^2/s$. 
As in K\"{a}mpf \cite{Kampf}, for the mesh we use the grid spacing of $\Delta x=5 m, \Delta z=2 m,$ and a time step $\Delta=0.1 s$. 

\begin{figure}[h]
\centering
$HDG_4$\hspace{6cm}K\"{a}mpf\\
10 minutes\\
\includegraphics[scale=0.2]{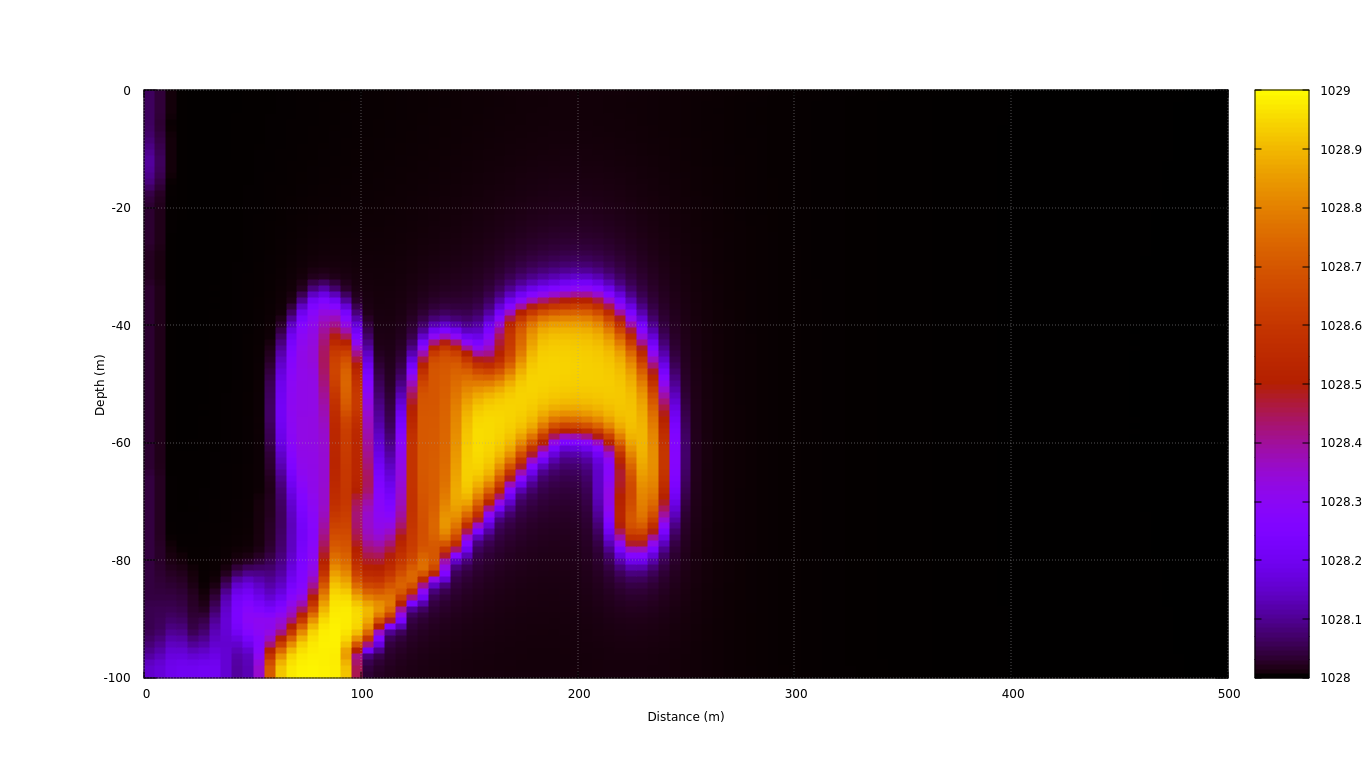}
\includegraphics[scale=0.2]{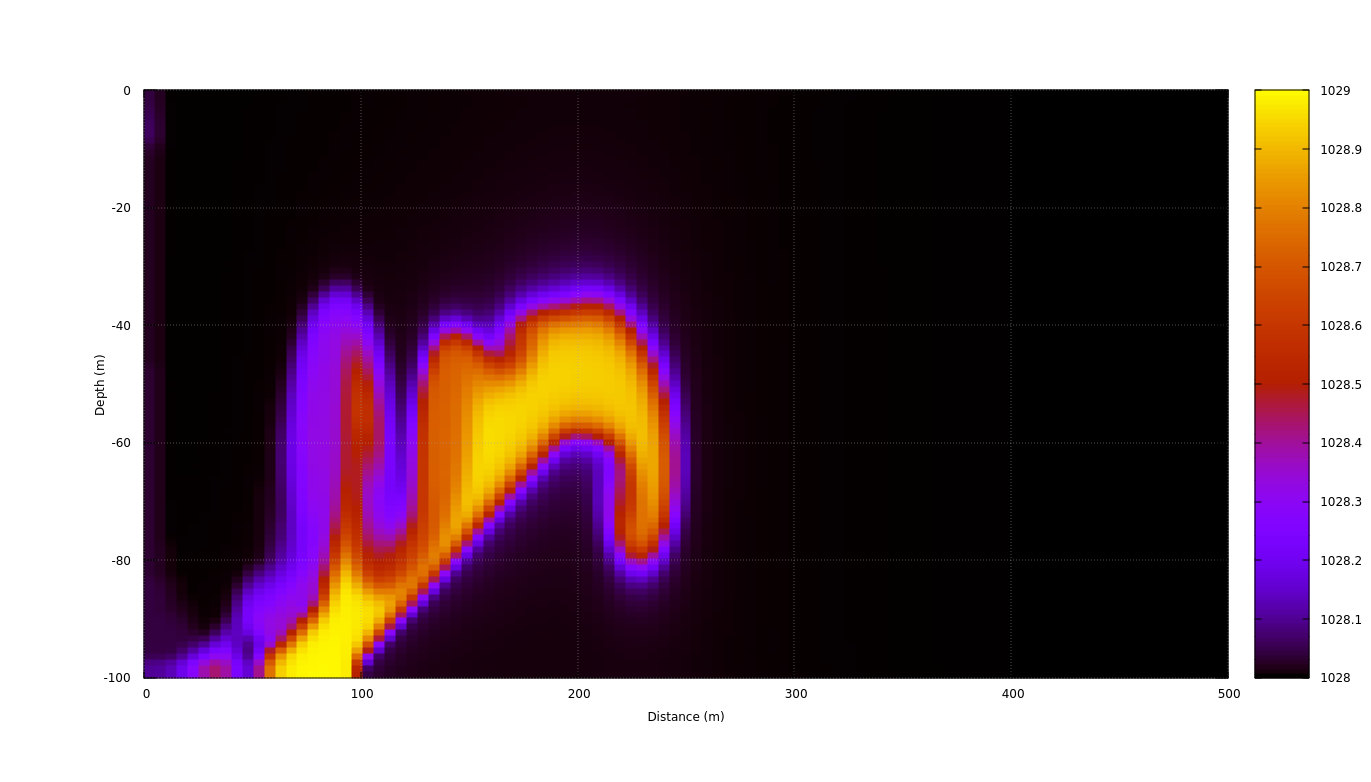}\\
15 minutes\\
\includegraphics[scale=0.2]{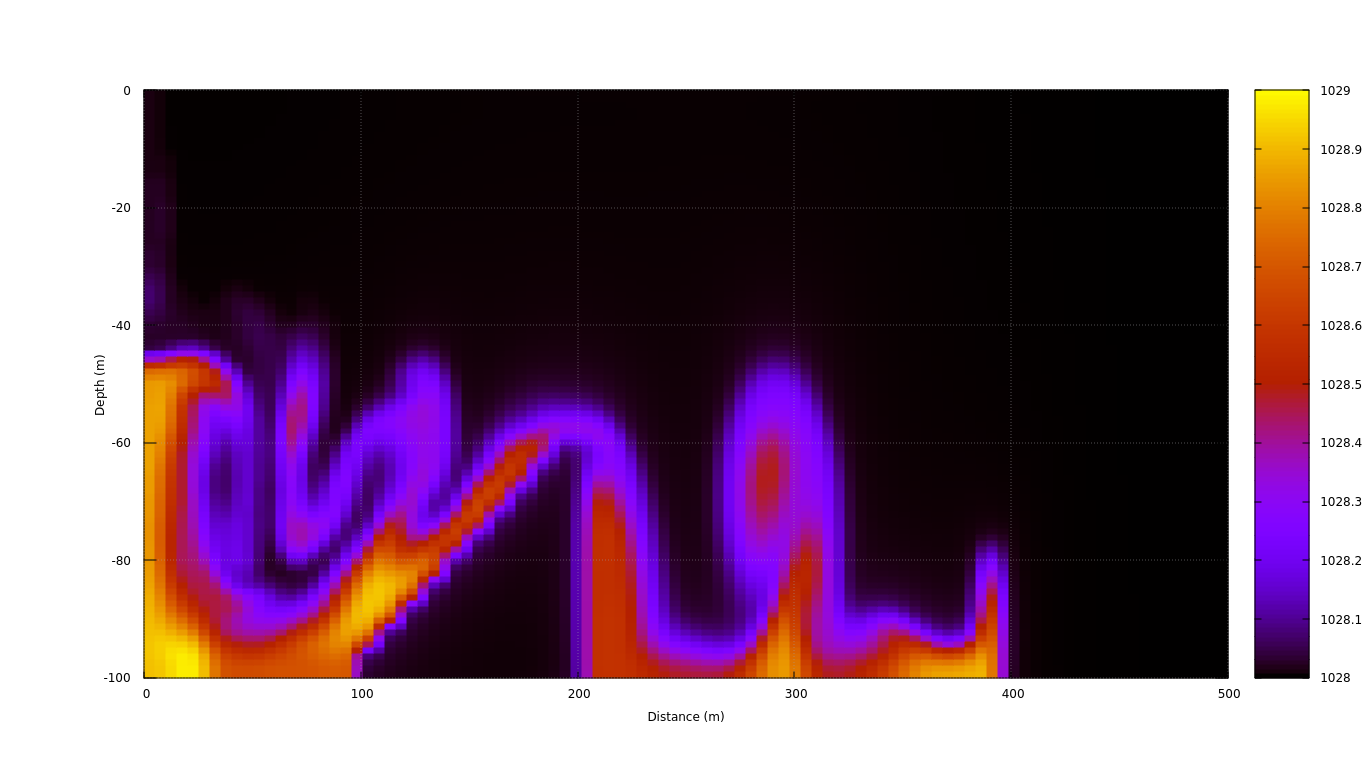}
\includegraphics[scale=0.2]{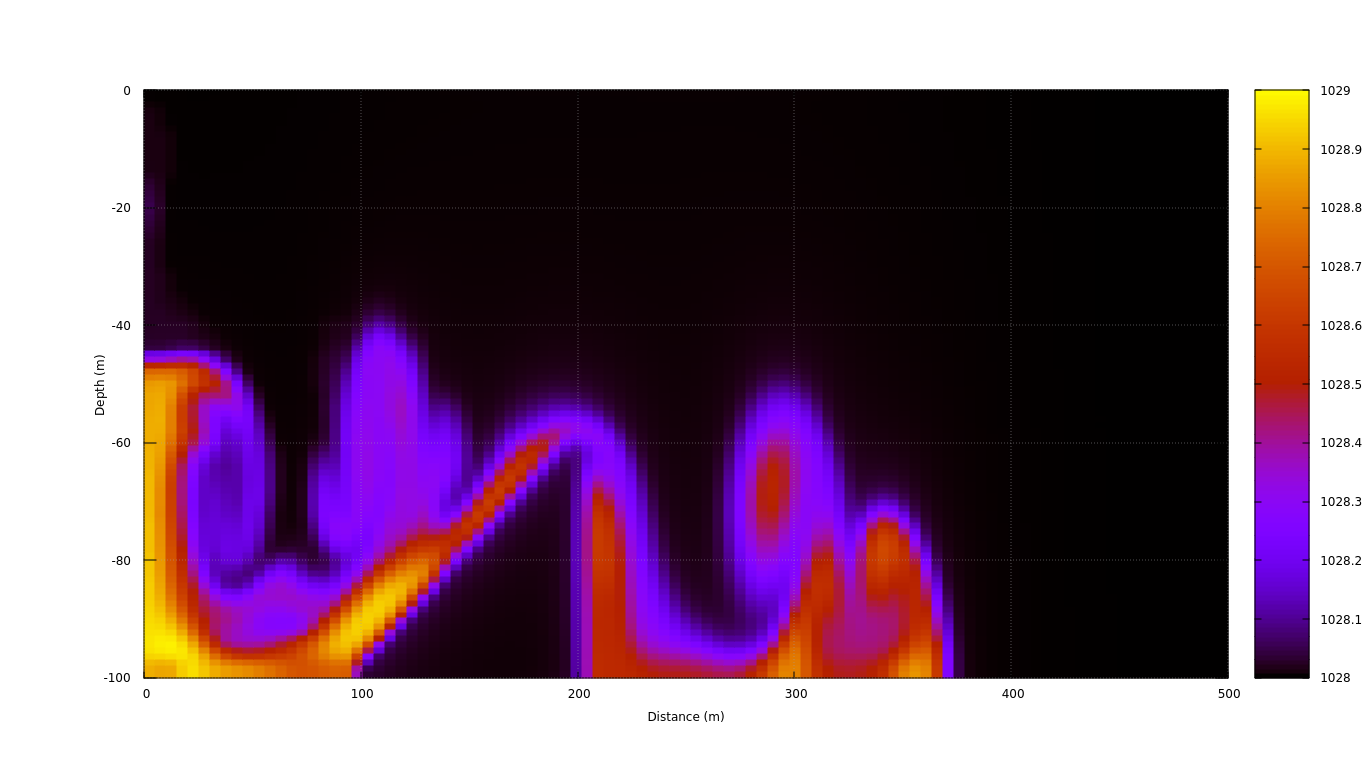}\\
20 minutes\\
\includegraphics[scale=0.2]{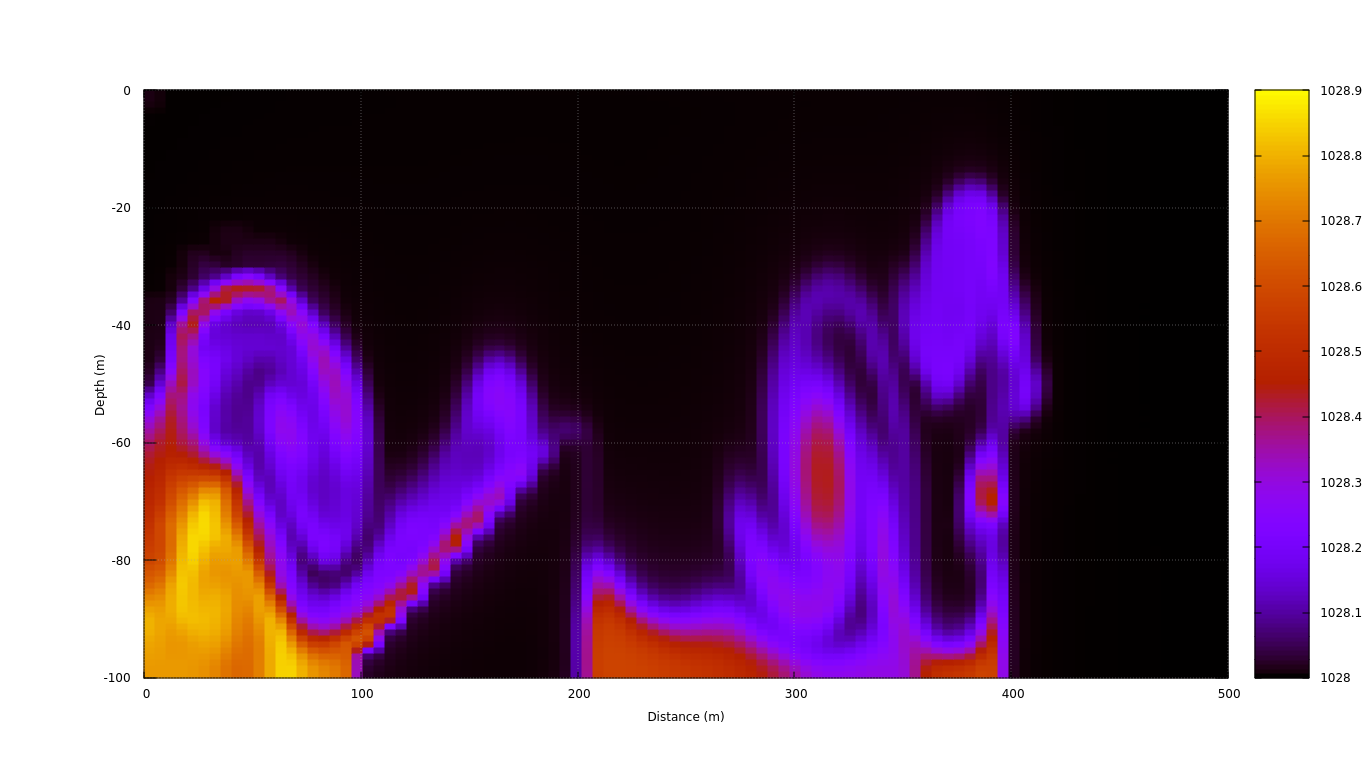}
\includegraphics[scale=0.2]{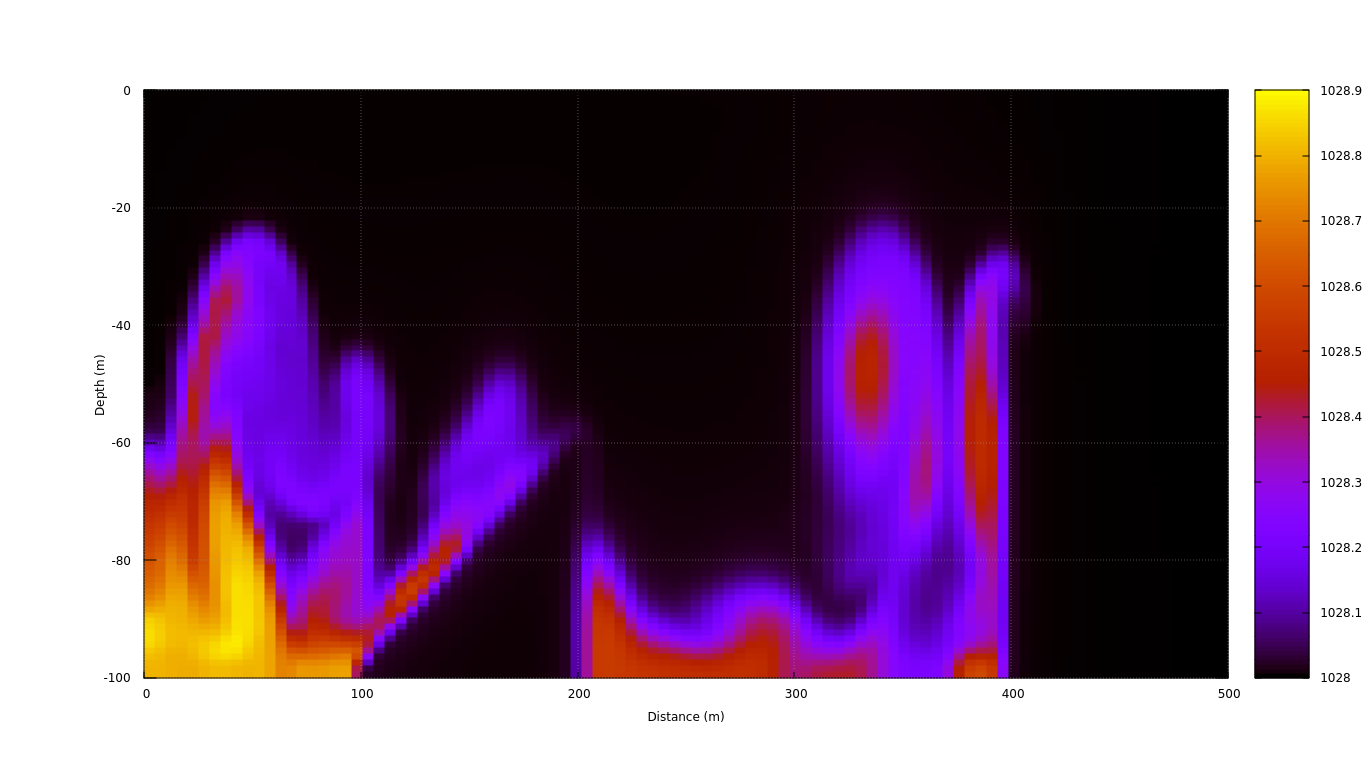}
\caption{Comparison of the numerical solution of $HDG_4$ and the FD2 algorithm for the variable density problem.}
\label{densidad_mueve_b}
\end{figure}

Numerical solutions are presented in the figure (\ref{densidad_mueve_b}). Therein, a comparison of the solution obtained with the $HDG_4$ scheme is compared with the second order finite difference method (FD2) in in K\"{a}mpf \cite{Kampf}.

\begin{figure}[h]
\centering
$HDG_2$\hspace{6cm}K\"{a}mpf\\
10 minutes\\
\includegraphics[scale=0.2]{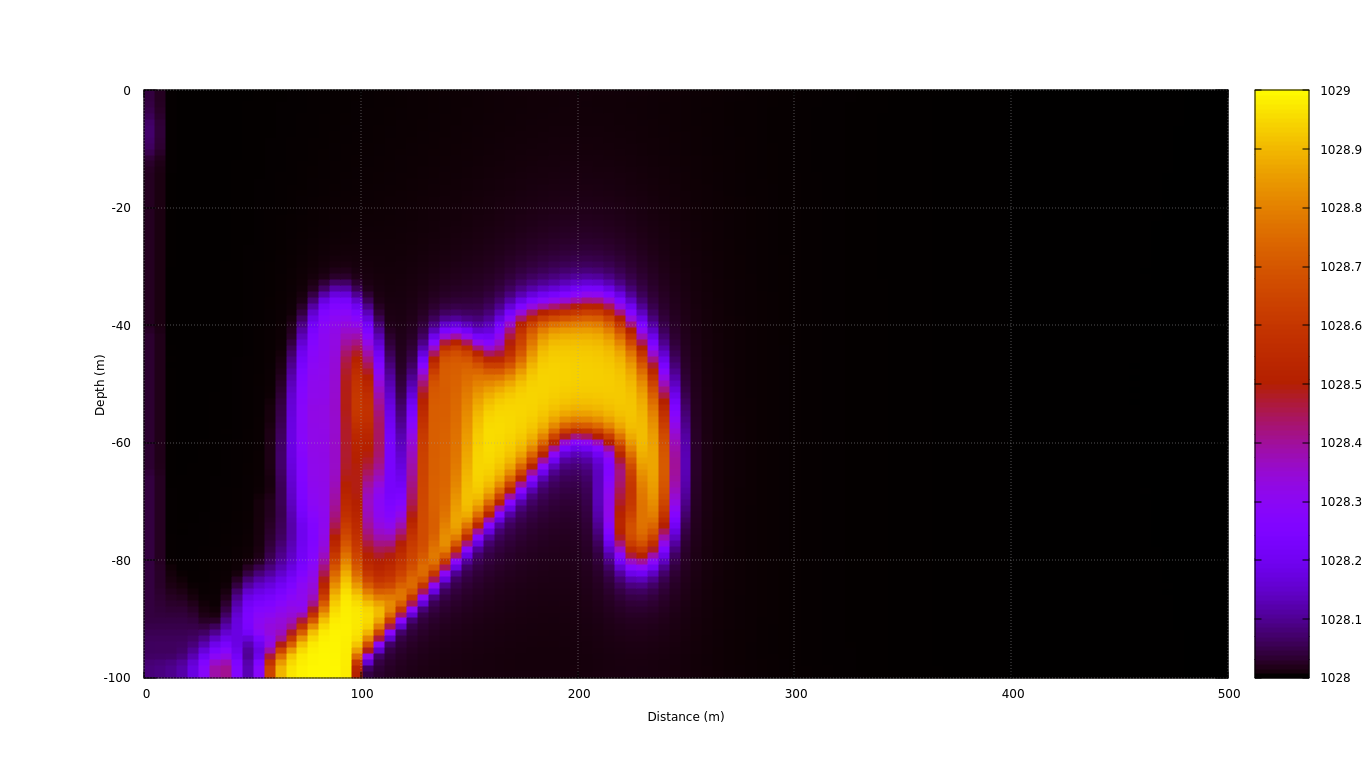}
\includegraphics[scale=0.2]{krho2v10.png}
15 minutes\\
\includegraphics[scale=0.2]{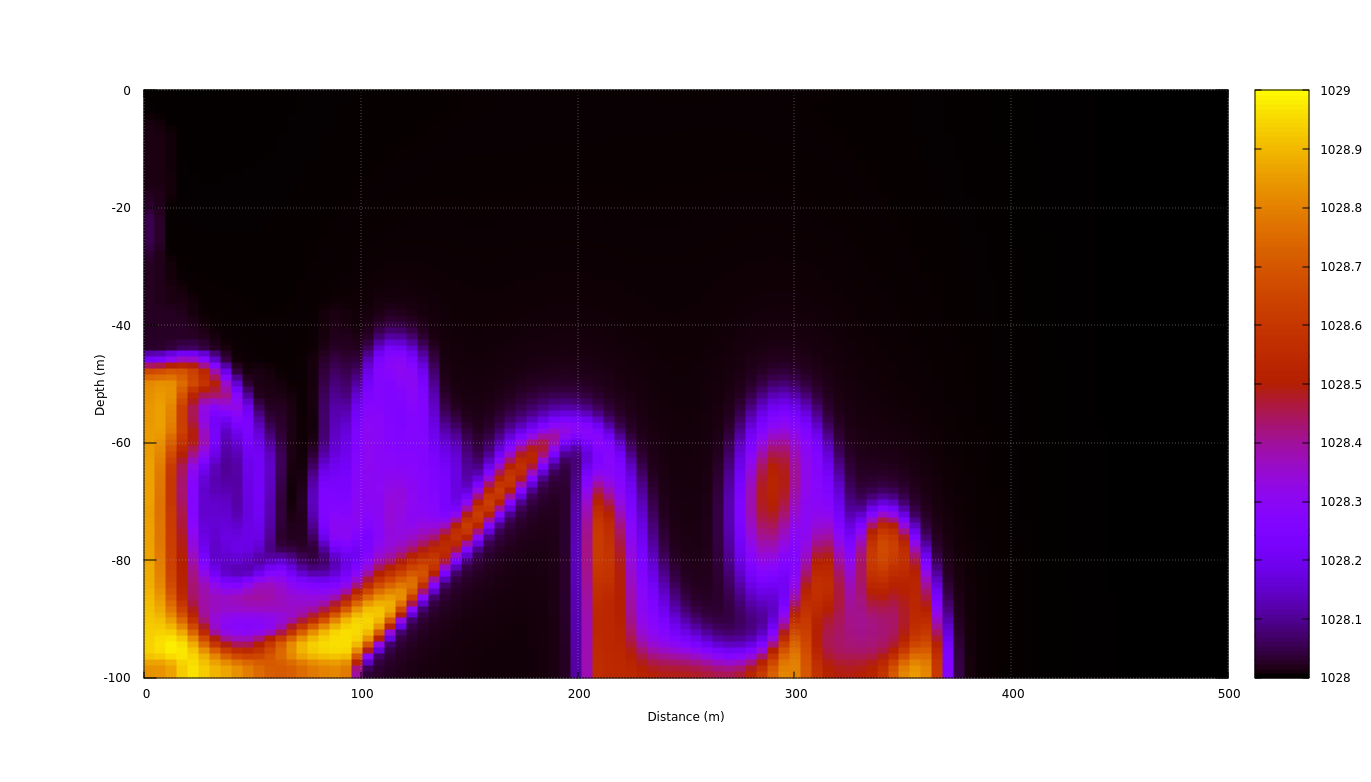}
\includegraphics[scale=0.2]{krho2v15.png}\\
20 minutes\\
\includegraphics[scale=0.2]{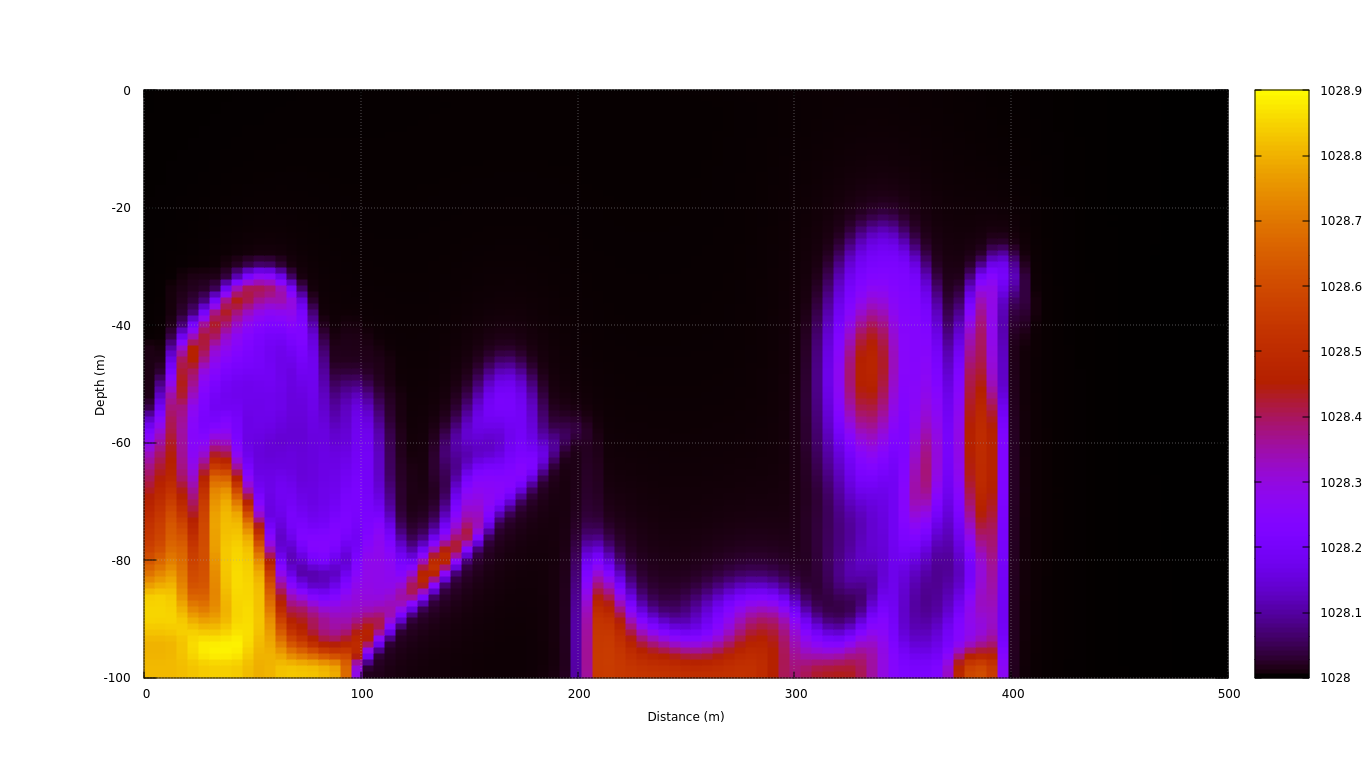}
\includegraphics[scale=0.2]{krho2v20.png}
\caption{Comparison of the numerical solution of $HDG_2$ and the K\"{a}mpf's algorithm for the variable density problem. }
\label{densidad_mueve_b_dg2}
\end{figure}

We also solve the problem with a second order Discontinuous Galerkin method, $HDG_2$. The approximation is comparable to FD2, this is illustrated in table \ref{tabla_hdg2_kampf_mvb}, using the $HDG_4$ as the true solution. The graphical comparison in figure \ref{densidad_mueve_b_dg2}.\\

\begin{table}[h]
\centering
\begin{tabular}{|c|c|c|}
\hline
Time (min)&$\mathcal{E}_2{(HDG_2)}$&$\mathcal{E}_2$(K\"{a}mpf's)\\
\hline
7&$6.55\times 10^{-5}$&$6.77\times 10^{-5}$\\
\hline
9&$0.00010759$&$0.000107866$\\
\hline
13&$0.000104016$&$0.000107906$\\
\hline
\end{tabular}
\caption{Comparison of the $HDG_2$ and K\"{a}mpf algorithms for the variable density problem  taking the $HDG_4$ as reference.}
\label{tabla_hdg2_kampf_mvb}
\end{table}

\section{Conclusions}

We have proposed a Hybrid Discontinuous Galerkin method for solutions of elliptic problems. It reduces the number of coupled unknowns arising from the classical Discontinuous Galerkin Method.  Moreover, it does not require penalization terms as typical DG solutions of elliptic problems. A high order implementation provides fast and robust solutions.

These features are illustrated on elliptic problems arising from Coastal Ocean Modeling. Even in the sector channel problem involving an
advection term, results are highly satisfactory.

In a pressure projection method for unsteady Navier-Stokes equations, the HDG elliptic solver  can be used to solve the Poisson problems in each time step. Thus accelerating computations. As a primer we have considered two slice non hydrostatic problems from the literature. It is shown that straightforward  modifications for the pressure correction method, yield accurate solutions.

We stress that the simplest grids and time advancing schemes have been used. The extension to more sophisticated situations is straightforward.

Of current and future interest, is to  apply this HDG approach to more sophisticated unsteady problems. We shall report on this elsewhere.

\end{document}